\documentclass[leqno,a4paper,10pt]{amsart}
\usepackage{amsmath,amsthm,amssymb,amsfonts,enumerate,color}

\date{September 28, 2010 - Revised version: February 3, 2011}

\newcommand{\norm}[1]{\left\Vert#1\right\Vert}
\newcommand{\abs}[1]{\left\vert#1\right\vert}
\newcommand{\set}[1]{\left\{#1\right\}}
\newcommand{\Real}{\mathbb{R}}
\newcommand{\card}{\operatorname{card}}
\newcommand{\T}{\mathcal{T}}
\renewcommand{\L}{\mathcal{L}}
\renewcommand{\P}{\mathcal{P}}

\newcommand{\supp}{\operatorname{supp}}

\newtheorem{thm}{Theorem}[section]

\newtheorem{lem}[thm]{Lemma}

\theoremstyle{definition}

\newtheorem{rem}[thm]{Remark}

\numberwithin{equation}{section}

\author[I. Abu-Falahah]{Ibraheem Abu-Falahah}
\address{Departamento de Matem\'aticas \\
          Facultad de Ciencias \\
          Universidad Au\-t\'o\-no\-ma de Madrid \\
          28049 Madrid, Spain}
\email{ibraheem.abufalahah@uam.es}

\author[P. R. Stinga]{Pablo Ra\'ul Stinga}
\address{Departamento de Matem\'aticas \\
          Facultad de Ciencias \\
          Universidad Au\-t\'o\-no\-ma de Madrid \\
          28049 Madrid, Spain}
\email{pablo.stinga@uam.es}

\author[J. L. Torrea]{Jos\'e L. Torrea}
\address{Departamento de Matem\'aticas \\
          Facultad de Ciencias \\
          Universidad Au\-t\'o\-no\-ma de Madrid \\
          and ICMAT-CSIC-UAM-UCM-UC3M \\
          28049 Madrid, Spain}
\email{joseluis.torrea@uam.es}

\thanks{Research supported by Ministerio de Ciencia e Innovaci\'{o}n de Espa\~{n}a MTM2008-06621-C02-01}



\begin{document}

\title{Square functions associated to Schr\"odinger operators}

\begin{abstract}
We characterize geometric properties of Banach spaces in terms of boundedness of square functions associated to general Schr\"odinger operators of the form $\L=-\Delta+V$, where the nonnegative potential $V$ satisfies a reverse H\"older inequality.  The main idea is to sharpen the well known localization method introduced  by  Z. Shen.  Our results can be regarded as alternative  proofs of  the boundedness in  $H^1$, $L^p$ and $BMO$  of classical $\L$-square functions.
\end{abstract}

\maketitle

\section{Introduction}

Consider the time independent Schr\"odinger operator in $\Real^d$, $d\geq3$,
\begin{equation}\label{L}
\L:=-\Delta+V,
\end{equation}
where the nonnegative potential $V$ satisfies a reverse H\"older inequality for some $s>d/2$, see \eqref{RHs}.

Let $X$ be a Banach space and let $\set{\P_t}_{t>0}=\{e^{-t\sqrt{\L}}\}_{t>0}$ be the (subordinated) Poisson semigroup   associated to $\L$, see \eqref{subordinacion}. For $2\leq q<\infty$ consider the generalized square function
\begin{equation}\label{defig}
g^{\L,q}f(x)=\left(\int_0^\infty\norm{t\frac{\partial\P_tf(x)}{\partial t}}_X^q~\frac{dt}{t}\right)^{1/q}=\norm{t\partial_t\P_tf(x)}_{L^q_X((0,\infty),\frac{dt}{t})},\qquad x\in\Real^d.
\end{equation}
By using the method described below we prove the following Theorem.

\

\noindent \textbf{Theorem A.}
\textit{Let $X$ be a Banach space and $2\leq q<\infty$. The following statements are equivalent.
\begin{enumerate}[(i)]
    \item $X$ admits an equivalent norm for which it is $q$-uniformly convex.
    \item The operator $g^{\L,q}$ maps $BMO_{\L,X}$ into $BMO_{\L}$.
    \item The operator $g^{\L,q}$ maps $L^p_X(\Real^d)$ into $L^p(\Real^d)$, for any $p$ in the range $1<p<\infty$.
    \item The operator $g^{\L,q}$ maps $L^1_X(\Real^d)$ into weak-$L^1(\Real^d)$.
    \item The operator $g^{\L,q}$ maps $H^1_{\L,X}$ into $L^1(\Real^d)$.
    \item For every $f\in L^1_X(\Real^d)$, $g^{\L,q}f(x)<\infty$ for almost every $x\in\Real^d$.
\end{enumerate}}

\

In 1995 Z. Shen proved $L^p$-boundedness of the Riesz transforms associated to the operator $\L$, see \cite{Shen}. The main idea in that paper is to break the kernels of the operators into ``local''  and ``global'' parts (close to the diagonal and far from the diagonal according to a certain distance $\rho(x)$ related to $\L$). Such a paper, a nice and exhaustive piece of mathematics, has became a classic and it has been a source of inspiration for a lot of manuscripts regarding Harmonic Analysis of operators associated to \eqref{L}. However, when these operators are defined with some formula involving the heat semigroup (as in the case of the maximal operator $\sup_{t>0}\abs{e^{-t\L}f}$ and the square function $\big(\int_0^\infty\abs{t\partial_te^{-t\L}f}^2\tfrac{dt}{t}\big)^{1/2}$) the word ``locally'' usually refers to the parameter $t$ of $e^{-t\L}$ being small and controlled in some sense by $\rho(x)$, see Dziuba\'nski et al., \cite{DGMTZ}.

Beyond the characterization of $q$-uniformly convex Banach spaces through boundedness properties of $\L$-square functions, we have another purpose. Namely enlighten the ``localization'' technique by sharpening the method introduced in \cite{Shen} in order to avoid the manipulations with the parameter $t$. At the same time, we get a unified approach to prove $H^1$, $L^p$ and $BMO$ boundedness results for classical Harmonic Analysis operators associated to $\L$. Observe that, in particular, Theorem A gives an alternative proof of the boundedness of $g^{\L,2}$ in the scalar case.

Let us briefly describe the procedure that within the paper is developed in detail for the case of the square function $g^{\L,q}$ acting on vector valued functions.

\

\textbf{Description of the method.} \textit{Let $\rho(x)$ be the auxiliary critical radii function determined by the potential $V$, see \eqref{critical}, and $N$ be the region consisting of points $(x,y)\in\Real^d\times\Real^d$ such that $\abs{x-y}\leq\rho(x)$. Given a linear operator associated to $\L$, that we denote by $T^\L$, let $T^\Delta$ be the parallel operator associated to the classical Laplacian $-\Delta$. Define the localized operator $T^\L_{\mathrm{loc}}f(x):=T^\L\left(\chi_N(x,\cdot)f(\cdot)\right)(x)$ and analogously $T^\Delta_{\mathrm{loc}}f(x)$. Then $T_{\mathrm{loc}}^\Delta$ inherits the $L^p$-boundedness properties from the operator $T^\Delta$. Even more, if $T_{\mathrm{loc}}^\Delta$ is bounded in $L^p$ then it is also bounded in $BMO_\L$. In other words, the operator $T_{\mathrm{loc}}^\Delta$ behaves as a natural operator associated to $\L$. Now the method finishes by observing that the difference operators $T_{\mathrm{loc}}^\Delta-T_{\mathrm{loc}}^{\L}$ and $T_{\mathrm{loc}}^{\L}-T^{\L}$ are bounded from $L^p$ into $L^p$ for $1\leq p\leq\infty$ and from $BMO_\L$ into $L^\infty$.}

\

In order to unify the method we consider a ``local'' part, defined through $\rho$, where the cutting acts on the heat kernel. This idea allows us to handle any operator defined via a formula involving the heat kernel, as for example Riesz transforms, square functions, etc.

Besides the unification, we believe that our main contribution is to show that the local part (given in terms of the distance $\rho$) of an operator associated to the standard Laplacian $-\Delta$ shares the natural boundedness properties with the corresponding operators associated to $\L$. We must emphasize how surprising this phenomenon is in the case of boundedness in $BMO$. See Theorems \ref{herencia} and \ref{trasplante} in Section \ref{technical} for the case of the $g$-function. The general ideas are summarized in Remarks \ref{metodo1} and \ref{metodo2}. Observe that the localized operator has always a rough kernel, see Remark \ref{metodo2}, and it is not clear a priori how to prove the necessary smoothness properties in order to get the desired boundedness in $BMO$.

Here $L^p_X(\Real^d)$ denotes the usual $L^p$-space of Bochner-Lebesgue $p$-integrable functions on $\Real^d$ with values in $X$. The spaces $H^1_{\L,X}$ and $BMO_{\L,X}$ are defined in the same way as in the scalar case just by replacing the absolute value of $\mathbb{C}$ by the norm of $X$, see \eqref{h1} and \eqref{bmo}. For the definition of $q$-uniform convexity we refer to Section \ref{Section:Proof}. Throughout the paper the letter $C$ denotes a positive constant that may change in each appearance and does not depend on the significant quantities.

The paper is organized as follows. We collect in Section \ref{pre} the preliminary results already known in the context of Schr\"odinger operators. Section \ref{technical} contains the technical results needed for the application of the method. Finally Section \ref{Section:Proof} is devoted to the proof of Theorem A.

\vskip0.5cm
\noindent\textbf{Acknowledgement.} We are very grateful to the referee for the thorough revision of the original manuscript. His exhaustive and detailed comments, together with his helpful suggestions, certainly helped us to strongly improve the presentation of the paper.

\section{Preliminaries}\label{pre}

The nonnegative potential $V$ in \eqref{L} satisfies a reverse H\"older inequality for some $s>d/2$; that is, there exists a constant $C=C(s,V)$ such that
\begin{equation}\label{RHs}
\left(\frac{1}{\abs{B}}\int_BV(y)^s~dy\right)^{1/s}\leq\frac{C}{\abs{B}}\int_BV(y)~dy,
\end{equation}
for all balls $B\subset\Real^d$. Associated to this potential, Shen defines in \cite{Shen} the critical radii function as
\begin{equation}\label{critical}
\rho(x):=\sup\Big\{r>0:\frac{1}{r^{d-2}}\int_{B(x,r)}V(y)~dy\leq1\Big\},\qquad x\in\Real^d.
\end{equation}
Some properties of this function $\rho$ are well known. We are particularly interested in the following.

\begin{lem}[see Lemma~1.4 in \cite{Shen}]\label{Lem:equiv rho}
There exist $c>0$ and $k_0\geq1$ so that for all $x,y\in\Real^d$
\begin{equation}\label{2}
c^{-1}\rho(x)\left(1+\frac{\abs{x-y}}{\rho(x)}\right)^{-k_0}\leq\rho(y)\leq c\rho(x)\left(1+\frac{\abs{x-y}}{\rho(x)}\right)^{\frac{k_0}{k_0+1}}.
\end{equation}
In particular, there exists a positive constant $C_1<1$ such that
$$\hbox{if}\quad\abs{x-y}\leq\rho(x)\quad\hbox{then}\quad C_1\rho(x)<\rho(y)<C_1^{-1}\rho(x).$$
\end{lem}

\begin{lem}[see Lemma~2.3 in \cite{Dziubanski-Zienkiewicz}]\label{bolas}
There exists a sequence of points $\set{x_k}_{k=1}^\infty$ in $\Real^d$ such that the family of balls $\set{Q_k}_{k=1}^\infty$ defined by $Q_k:= B(x_k,\rho(x_k))$ satisfy
\begin{itemize}
    \item $\bigcup_k Q_k=\Real^d$;
    \item There exists $N=N(\rho)$ so that, for every $k\geq1$, $\card\set{j:2Q_j\cap2Q_k\neq\emptyset}\leq N$;
\end{itemize}
where for a ball $B$ and a positive number $c$ we denote by $cB$ the ball with the same center as $B$ and radius $c$ times the radius of $B$.
\end{lem}

Let $\set{\T_t}_{t>0}$ be the heat--diffusion semigroup associated to $\L$ acting on $X$-valued functions:
\begin{equation}\label{heatL}
\T_tf(x)\equiv e^{-t\L}f(x)=\int_{\Real^d}k_t(x,y)f(y)~dy,\qquad f\in L^2_X(\Real^d),~x\in\Real^d,~t>0.
\end{equation}
The following Lemmas are known.

\begin{lem}[see \cite{Dziubanski-Zienkiewicz Hp,Kurata}]\label{Lem:cota heat L}
For every $\alpha>0$ there exists a constant $C_\alpha$ such that
\begin{equation}\label{1}
0\leq k_t(x,y)\leq C_\alpha\frac{1}{t^{d/2}}~e^{-\frac{\abs{x-y}^2}{5t}}\left(1+\frac{\sqrt{t}}{\rho(x)}+\frac{\sqrt{t}}{\rho(y)}\right)^{-\alpha},
\end{equation}
for all $x,y\in\Real^d$, $t>0$.
\end{lem}

Let
$$h_t(x):=\frac{1}{(4\pi t)^{d/2}}~e^{-\frac{\abs{x}^2}{4t}},\qquad x\in\Real^d,~t>0,$$
be the kernel of the classical heat semigroup $\set{T_t}_{t>0}=\{e^{t\Delta}\}_{t>0}$ in $\Real^d$.

\begin{lem}[see Proposition 2.16 in \cite{Dziubanski-Zienkiewicz Hp}]\label{Lem:Schwartz}
There exists a nonnegative Schwartz class function $\omega$ in $\Real^d$ such that
\begin{equation}\label{3}
\abs{k_t(x,y)-h_t(x-y)}\leq\left(\frac{\sqrt{t}}{\rho(x)}\right)^\delta\omega_t(x-y),\qquad x,y\in\Real^d,~t>0,
\end{equation}
where $\omega_t(x-y):=t^{-d/2}\omega\left((x-y)/\sqrt{t}\right)$ and
\begin{equation}\label{delta}
\delta:=2-\frac{d}{s}>0.
\end{equation}
\end{lem}

Given the heat semigroup \eqref{heatL}, the Poisson semigroup associated to $\L$ is obtained through Bochner's subordination formula, see \cite{SteinTopics}:
\begin{equation}\label{subordinacion}
\P_tf(x)\equiv e^{-t\sqrt{\L}}f(x)=\frac{t}{2\sqrt{\pi}}\int_0^\infty\frac{e^{-t^2/(4u)}}{u^{3/2}}~\T_uf(x)~du,\qquad x\in\Real^d,~t>0.
\end{equation}
With this we define, for $2\leq q<\infty$, the square function related to $\L$ as in \eqref{defig}.

\begin{rem}[Notational convention]\label{notacion}
The Poisson semigroup associated to the classical Laplace operator in $\Real^d$ will be denoted by $\set{P_t}_{t>0}=\{e^{-t\sqrt{-\Delta}}\}_{t>0}$. Recall that $P_tf(x)=P_t\ast f(x)$, where
$$P_t(x)=c_d\frac{t}{(t^2+\abs{x}^2)^{\frac{d+1}{2}}},\qquad x\in\Real^d,~t>0.$$
The square function considered in \eqref{defig} will be denoted by $g^{\Delta,q}f$ when replacing $\P_tf$ by $P_tf$.
\end{rem}

A locally integrable function $f:\Real^d\to X$ is in $BMO_{\L,X}$ whenever there exists a constant $C$ such that
\begin{enumerate}
\item[(i)] $\displaystyle\frac{1}{\abs{B}}\int_B\norm{f(x)-f_B}_Xdx\leq C$, for every ball $B$ in $\Real^d$, and
\item[(ii)] $\displaystyle\frac{1}{\abs{B}}\int_B\norm{f(x)}_Xdx\leq C$, for every $B=B(x_0,r_0)$, where $x_0\in\Real^d$ and $r_0\geq\rho(x_0)$.
\end{enumerate}
As usual, $f_{B}:=\displaystyle\frac{1}{|B|}\int_Bf(x)~dx$, for every ball $B$ in $\Real^d$. The norm $\norm{f}_{BMO_{\L,X}}$ of $f$ is defined as
\begin{eqnarray}\label{bmo}
\norm{f}_{BMO_{\L,X}}=\inf\set{C\geq0:\hbox{(i)~and~(ii)~hold}}.
\end{eqnarray}
Let us note that if $\mathrm{(ii)}$ is true for some ball $B$ then $\mathrm{(i)}$ holds true for the same ball, so we might ask to $\mathrm{(i)}$ only for balls with radii smaller than $\rho(x_0)$. By using the classical John-Nirenberg inequality it can be seen that if in (i) and (ii) $L^1_X$-norms are replaced by $L^p_X$-norms, for $1<p<\infty$, then the space $BMO_{\L,X}$ does not change and equivalent norms appear, see \cite[Corollary~3]{DGMTZ}.

We define the vector-valued atomic Hardy space related to $\L$ following the scalar-valued definition in \cite{Dziubanski-Zienkiewicz}. A function $a:\Real^d\to X$ is an $H^1_{\L,X}$-atom associated with a ball $B(x_0,r)$ when $\supp a\subset B(x_0,r)$,
\begin{equation}\label{prop atom}
\norm{a}_{L^\infty_X(\Real^d)}\leq\frac{1}{\abs{B(x_0,r)}},
\end{equation}
and, in addition,
\begin{equation}\label{prop atom media}
\int_{\Real^d}a(x)~dx=0,\qquad\hbox{whenever}~0<r<\rho(x_0).
\end{equation}
An $X$-valued integrable function $f$ in $\Real^d$ belongs to $H^1_{\L,X}$ if and only if it can be written as $f=\sum_j\lambda_ja_j$, where $a_j$ are $H^1_{\L,X}$-atoms and $\sum_j\abs{\lambda_j}<\infty$. The norm is given by
\begin{equation}\label{h1}
\norm{f}_{H^1_{\L,X}}=\inf\Big\{\sum_j\abs{\lambda_j}:f=\sum_j\lambda_ja_j\Big\}.
\end{equation}
In \cite{DGMTZ} it is shown that $BMO_{\L}$ is the dual space of $H^1_{\L}$.

\section{Technical Lemmas}\label{technical}

As we said in the description of our method, the following region $N$ will play a fundamental role:
$$N:=\set{(x,y)\in\Real^d\times\Real^d:\abs{x-y}\leq\rho(x)}.$$
Given $N$ we define the ``global'' and ``local'' parts of the square function defined in (\ref{defig}) as
\begin{align}\label{global}
    g^{\L,q}_{\mathrm{glob}}f(x) &= g^{\L,q}\left(\chi_{N^c}(x,\cdot)f(\cdot)\right)(x)\qquad\hbox{and} \\
    g^{\L,q}_{\mathrm{loc}}f(x) &= g^{\L,q}f(x)-g^{\L,q}_{\mathrm{glob}}f(x)\nonumber.
\end{align}
Note that
\begin{equation}\label{eq:local control}
g^{\L,q}_{\mathrm{loc}}f(x)\leq g^{\L,q}(\chi_N(x,\cdot)f(\cdot))(x)\leq g^{\L,q}f(x)+g^{\L,q}_{\mathrm{glob}}f(x),\qquad\hbox{a.e.}~x\in\Real^d,
\end{equation}
or equivalently,
\begin{equation}\label{eq:dif con corte}
\abs{g^{\L,q}f(x)-g^{\L,q}(\chi_N(x,\cdot)f(\cdot))(x)}\leq g^{\L,q}_{\mathrm{glob}}f(x),\qquad\hbox{a.e.}~x\in\Real^d.
\end{equation}

\begin{lem}\label{lemglob}
Let $X$ be any Banach space and $\alpha>0$. Then for any $f\in\bigcup_{1\leq p\leq\infty}L^p_X(\Real^d)$ we have
$$g^{\L,q}_{\mathrm{glob}}f(x)\leq C\int_{\Real^d}L(x,y)\chi_{N^c}(x,y)\norm{f(y)}_Xdy,\qquad x\in\Real^d,$$
where $\displaystyle L(x,y)=\frac{\rho(x)^\alpha}{\abs{x-y}^{d+\alpha}}$, $x,y\in\Real^d$.
\end{lem}

\begin{proof}
Using Bochner's subordination formula \eqref{subordinacion} it can be checked that for any function $h$,
$$\norm{\partial_t\P_th(x)}_X\leq C\int_0^\infty\frac{e^{-t^2/(8u)}}{u^{3/2}}\norm{\T_uh(x)}_Xdu,$$
where we applied the inequality $r^\eta e^{-r}\leq C_\eta e^{-r/2}$, valid for $\eta\geq0$, $r>0$. Hence, by Minkowski's inequality,
\begin{align*}
    g^{\L,q}_{\mathrm{glob}}f(x) &\leq C\int_0^\infty\Big\|t\frac{e^{-t^2/(8u)}}{u^{3/2}}\Big\|_{L^q((0,\infty),\frac{dt}{t})} \norm{\T_u(\chi_{N^c}(x,\cdot)f(\cdot))(x)}_Xdu \\
    &= C\int_0^\infty\norm{\T_u(\chi_{N^c}(x,\cdot)f(\cdot))(x)}_X\frac{du}{u}\leq C\int_0^\infty\int_{\Real^d}k_u(x,y)\chi_{N^c}(x,y)\norm{f(y)}_Xdy~\frac{du}{u}.
\end{align*}
From \eqref{1} of Lemma \ref{Lem:cota heat L} and the change of variables $r=\frac{\abs{x-y}^2}{cu}$ we get
$$\int_0^\infty k_u(x,y)~\frac{du}{u}\leq C\int_0^\infty\frac{1}{u^{d/2}}~e^{-\frac{\abs{x-y}^2}{cu}}\left(\frac{\rho(x)}{\sqrt{u}}\right)^\alpha\frac{du}{u} =C\frac{\rho(x)^\alpha}{\abs{x-y}^{d+\alpha}}\int_0^\infty r^{\frac{d+\alpha}{2}}e^{-r}\frac{dr}{r}.$$
\end{proof}

\begin{lem}\label{Lem:ltodop}
Let $X$ be any Banach space. Then the global operator $g^{\L,q}_{\mathrm{glob}}$ maps
\begin{itemize}
    \item[(a)] $L^p_X(\Real^d)$ into $L^p(\Real^d)$ for any $p$, $1\leq p\leq\infty$,
    \item[(b)] $BMO_{\L,X}$ into $L^\infty(\Real^d)$, and
    \item[(c)] $H^1_{\L,X}$ into $L^1(\Real^d)$.
\end{itemize}
\end{lem}

\begin{proof}
Let $L(x,y)$, $x,y\in\Real^d$, be as in Lemma \ref{lemglob}. Observe that
$$\int_{\Real^d}L(x,y)\chi_{N^c}(x,y)~dy=\rho(x)^{\alpha}\int_{\abs{x-y}>\rho(x)}\frac{1}{\abs{x-y}^{d+\alpha}}~dy=C,$$
for all $x\in\Real^d$. On the other hand, by Lemma \ref{Lem:equiv rho}, there exists a positive number $\varepsilon<1$ such that
\begin{equation}\label{estimateL}
L(x,y)\leq C\frac{\rho(y)^{\alpha}}{\abs{x-y}^{d+{\alpha}}}\left(1+\frac{\abs{x-y}}{\rho(y)}\right)^{\varepsilon {\alpha}}\leq C\left(\frac{\rho(y)^{\alpha}}{\abs{x-y}^{d+{\alpha}}}+\frac{\rho(y)^{(1-\varepsilon)\alpha}}{\abs{x-y}^{d+(1-\varepsilon)\alpha}}\right).
\end{equation}
Assume that $\abs{x-y}>\rho(x)$. Then we claim that $\abs{x-y}\geq C\rho(y)$ for some positive constant $C$ depending on the constants $c$ and $k_0$ that appear in Lemma 2.1. Indeed, by Lemma 2.1 and the fact that $\frac{\abs{x-y}}{\rho(x)}\geq1$ and $\frac{k_0}{k_0+1}\leq 1$, we have
$$\rho(y)\leq C\rho(x)\left(1+\left(\frac{\abs{x-y}}{\rho(x)}\right)^{\frac{k_0}{k_0+1}}\right)\leq C\rho(x)\left(1+\frac{\abs{x-y}}{\rho(x)}\right)\leq C\left(\rho(x)+\abs{x-y}\right)\leq 2C\abs{x-y}.$$
This together with \eqref{estimateL} give us $\displaystyle \int_{\Real^d}L(x,y)\chi_{N^c}(x,y)~dx\leq C$. Hence the operator given by the kernel $L(x,y)\chi_{N^c}(x,y)$ maps $L^p_X(\Real^d)$ into $L^p(\Real^d)$ for every $p$, $1\leq p\leq\infty$. Using Lemma \ref{lemglob} we get $\mathrm{(a)}$.

In order to see $\mathrm{(b)}$ we observe that for a function $f$ in $BMO_{\L,X}$, by Lemma \ref{lemglob},
\begin{align*}
    g^{\L,q}_{\mathrm{glob}}f(x) &\leq C\rho(x)^{\alpha}\sum_{j=0}^\infty\int_{2^j\rho(x)<\abs{x-y}\leq2^{j+1}\rho(x)}\frac{1}{\abs{x-y}^{d+\alpha}}\norm{f(y)}_Xdy \\
    &\leq C\rho(x)^\alpha\sum_{j=0}^\infty\frac{1}{(2^j\rho(x))^{d+\alpha}}\int_{\abs{x-y}\leq2^{j+1}\rho(x)}\norm{f(y)}_Xdy \\
    &= C\sum_{j=0}^\infty\frac{1}{2^{j\alpha}}\frac{1}{(2^{j+1}\rho(x))^d}\int_{\abs{x-y}\leq2^{j+1}\rho(x)}\norm{f(y)}_Xdy \\
    &\leq C\norm{f}_{BMO_{\L,X}}\sum_{j=0}^\infty\frac{1}{2^{j\alpha}}=C\norm{f}_{BMO_{\L,X}},\qquad\hbox{for all}~x\in\Real^d.
\end{align*}

For $\mathrm{(c)}$ just note that $H^1_{\L,X}\subset L^1_X(\Real^d)$ and then apply $\mathrm{(a)}$.
\end{proof}

\begin{lem}\label{lemloc}
Let $X$ be any Banach space. Then, for any strongly measurable $X$-valued function $f$,
$$\abs{g^{\L,q}_{\mathrm{loc}}f(x)-g^{\Delta,q}_{\mathrm{loc}}f(x)}\leq C\int_{\Real^d}M(x,y)\chi_N(x,y)\norm{f(y)}_Xdy,\qquad x\in\Real^d,$$
where $\displaystyle M(x,y)=\frac{\rho(x)^{-\delta}}{\abs{x-y}^{d-\delta}}$, for $x,y\in\Real^d$, and $\delta>0$ is given in \eqref{delta}.
\end{lem}

\begin{proof}
Proceeding as in the proof of Lemma \ref{lemglob} it is easy to check that
$$\abs{g^{\L,q}_{\mathrm{loc}}f(x)-g^{\Delta,q}_{\mathrm{loc}}f(x)}\leq C\int_0^\infty\int_{\Real^d}\abs{k_u(x,y)-h_u(x-y)}\chi_{N}(x,y)\norm{f(y)}_Xdy~\frac{du}{u}.$$
Using \eqref{3} in Lemma \ref{Lem:Schwartz} and the fact that $\omega$ is a rapidly decreasing function,
\begin{align*}
  \lefteqn{\int_0^\infty\abs{k_u(x,y)-h_u(x-y)}~\frac{du}{u} \leq C\rho(x)^{-\delta}\int_0^\infty \frac{1}{u^{(d-\delta)/2}}~\omega\left((x-y)/\sqrt{u}\right)~\frac{du}{u}} \\
  &\leq C\rho(x)^{-\delta}\left[\frac{1}{\abs{x-y}^{d-\delta+\varepsilon}} \int_0^{\abs{x-y}^2}\left(\frac{\abs{x-y}}{\sqrt{u}}\right)^{d-\delta+\varepsilon} \omega\left((x-y)/\sqrt{u}\right)~\frac{du}{u^{1-\varepsilon/2}} +\int_{\abs{x-y}^2}^\infty\frac{1}{u^{(d-\delta)/2}}~\frac{du}{u}\right] \\
  &\leq C\frac{\rho(x)^{-\delta}}{\abs{x-y}^{d-\delta}}\left[\frac{1}{\abs{x-y}^\varepsilon}\int_0^{\abs{x-y}^2}~ \frac{du}{u^{1-\varepsilon/2}}+1\right]=C\frac{\rho(x)^{-\delta}}{\abs{x-y}^{d-\delta}}.
\end{align*}
\end{proof}

\begin{lem}\label{difloc}
Let $X$ be any Banach space. Then the difference operator $g^{\L,q}_{\mathrm{loc}}-g^{\Delta,q}_{\mathrm{loc}}$ maps
\begin{itemize}
    \item[(a)] $L^p_X(\Real^d)$ into $L^p(\Real^d)$ for any $p$, $1\leq p\leq\infty$,
    \item[(b)] $BMO_{\L,X}$ into $L^\infty(\Real^d)$, and
    \item[(c)] $H^1_{\L,X}$ into $L^1(\Real^d)$.
\end{itemize}
\end{lem}

\begin{proof}
Let $M(x,y)$, $x,y\in\Real^d$, be as in Lemma \ref{lemloc}. First note that
$$\int_{\Real^d}M(x,y)\chi_N(x,y)~dy=\rho(x)^{-\delta}\int_{\abs{x-y}\leq\rho(x)}\frac{1}{\abs{x-y}^{d-\delta}}~dy=C,\qquad x\in\Real^d.$$
On the other hand, by Lemma \ref{Lem:equiv rho},
$$M(x,y)\leq \frac{C}{\abs{x-y}^{d-\delta}}~\rho(y)^{-\delta}\left(1+\frac{\abs{x-y}}{\rho(y)}\right)^{k_0\delta}\leq C\left(\frac{\rho(y)^{-\delta}}{\abs{x-y}^{d-\delta}}+ \frac{\rho(y)^{-(1+k_0)\delta}}{\abs{x-y}^{d-(1+k_0)\delta}}\right),$$
where $k_0\geq1$. This, and the fact that $\abs{x-y}>\rho(x)$ implies $\abs{x-y}>C\rho(y)$ (see the proof of Lemma \ref{Lem:ltodop}), give $\displaystyle\int_{\Real^d}M(x,y)\chi_N(x,y)~dx\leq C$ for all $y\in\Real^d$. Applying Lemma \ref{lemloc} we conclude $\mathrm{(a)}$ and as a consequence we also get $\mathrm{(c)}$.

We shall prove $\mathrm{(b)}$. Let $f\in BMO_{\L,X}$. Then
\begin{align*}
    \label{localBMO} \int_{\Real^d}M(x,y)\chi_N(x,y)\norm{f(y)}_Xdy &= \sum_{j=0}^\infty\int_{2^{-(j+1)}\rho(x)<\abs{x-y}\leq2^{-j}\rho(x)}\frac{\rho(x)^{-\delta}}{\abs{x-y}^{d-\delta}} \norm{f(y)}_Xdy \\
    &\leq C\sum_{j=0}^\infty\frac{1}{2^{j\delta}}\frac{1}{(2^{-j}\rho(x))^d}\int_{\abs{x-y}\leq2^{-j}\rho(x)}\norm{f(y)}_Xdy \\
    &\leq C\sum_{j=0}^\infty\frac{1}{2^{j\delta}}\left[\frac{1}{(2^{-j}\rho(x))^d}\int_{\abs{x-y}\leq2^{-j}\rho(x)}\norm{f(y)-f_{B(x,2^{-j}\rho(x))}}_X ~dy\right. \\
    &\qquad\quad \left. +\sum_{k=0}^{j-1}\left(\norm{f_{B(x,2^{-k}\rho(x))}-f_{B(x,2^{-(k+1)}\rho(x))}}_X\right)+\norm{f_{B(x,\rho(x))}}_X\right] \\
    &\leq C\sum_{j=0}^\infty\frac{1}{2^{j\delta}}\left[\norm{f}_{BMO_X}+j\norm{f}_{BMO_X}+\norm{f}_{BMO_{X,\L}}\right] \\
    &\leq C\norm{f}_{BMO_{\L,X}}\sum_{j=0}^\infty\frac{(j+2)}{2^{j\delta}}=C\norm{f}_{BMO_{\L,X}},
\end{align*}
for all $x\in\Real^d$. To finish use Lemma \ref{lemloc}.
\end{proof}

\begin{lem}\label{independiente}
Let $C_1$ be the constant that appears in Lemma \ref{Lem:equiv rho} and $\gamma>0$. Take $x,y\in\Real^d$ such that $\abs{x}<\gamma$ and $\abs{y}<\frac{C_1^2}{2}\rho(0)$. Then there exists a sufficiently large $R=R_\gamma>0$ for which
$\abs{\frac{x}{R}-y}<\rho\left(\frac{x}{R}\right)$.
\end{lem}

\begin{proof}
Lemma \ref{Lem:equiv rho} ensures that $C_1\rho(0)<\rho(y)<C_1^{-1}\rho(0)$. Let $R>0$ be such that $\abs{\frac{x}{R}-y}<C_1^2\rho(0)$ (it is enough to take $R>\frac{2\gamma}{C_1^2\rho(0)}$). Hence $\abs{\frac{x}{R}-y}<C_1\rho(y)<\rho(y)$. Once more using Lemma \ref{Lem:equiv rho} we obtain $\rho(y)<C_1^{-1}\rho\left(\frac{x}{R}\right)$ and therefore $\abs{\frac{x}{R}-y}<C_1C_1^{-1}\rho\left(\frac{x}{R}\right)=\rho\left(\frac{x}{R}\right)$.
\end{proof}

\begin{lem}\label{encoger}
Let $f$ be a function with compact support. For a real number $r$ denote by $f^r$ the dilation of $f$ defined by $f^r(x):=f(rx)$, $x\in\Real^d$. Then for any given $\gamma>0$ there exists $R>0$, depending on $\gamma$ and the support of $f$, such that
$$g^{\Delta,q}f(x)=g^{\Delta,q}\left(\chi_N(\tfrac{x}{R},\cdot)f^R(\cdot)\right)(\tfrac{x}{R}),\qquad\hbox{for all}~\abs{x}<\gamma.$$
\end{lem}

\begin{proof}
The scaling of the classical Poisson semigroup $P_tf^R(x/R)=P_{tR}f(x)$, $R>0$ (see Remark \ref{notacion}), implies that the square function satisfies $g^{\Delta,q}f(x)=g^{\Delta,q}f^R(x/R)$ for all $R>0$. In order to get the conclusion it is enough to take a sufficiently large $R$ such that the support of $f^R$ is contained in $B(0,\frac{C_1^2}{2}\rho(0))$ and such that Lemma \ref{independiente} can be applied.
\end{proof}

The following result establishes that the boundedness in $L^p$ of the square function $g^{\Delta,q}$ related to the Laplacian $-\Delta$ implies the same type of boundedness for the $\rho$-localized operator $g^{\Delta,q}_{\mathrm{loc}}$. In fact this is a fairly general property: see Remark \ref{metodo1} below.

\begin{thm}\label{herencia}
Assume that $g^{\Delta,q}$ maps $L^p_X(\Real^d)$ into $L^p(\Real^d)$ for some $p$, $1<p<\infty$ (resp. $L^1_X(\Real^d)$ into weak-$L^1(\Real^d)$). Then the operator $f\longmapsto t\partial_tP_t(\chi_{N}(x,\cdot)f(\cdot))(x)$, $x\in\Real^d$, $t>0$, maps $L^p_X(\Real^d)$ into $L^p_{L^q_X((0,\infty),\frac{dt}{t})}(\Real^d)$ (resp. $L^1_X(\Real^d)$ into weak-$L^1_{L^q_X((0,\infty),\frac{dt}{t})}(\Real^d)$). In particular $g^{\Delta,q}_{\mathrm{loc}}$ maps $L^p_X(\Real^d)$ into $L^p(\Real^d)$ (resp. $L^1_X(\Real^d)$ into weak-$L^1(\Real^d)$).

Moreover, if for every function $f\in L^1_X(\Real^d)$ we have $g^{\Delta,q}f(x)<\infty$ for almost all $x\in\Real^d$, then $\norm{t\partial_tP_t(\chi_N(x,\cdot)f(\cdot))(x)}_{L_X^q((0,\infty),\frac{dt}{t})}<\infty$ for almost all $x\in\Real^d$.
\end{thm}

\begin{proof}
We shall prove only the boundedness in $L^p$. We leave to the reader the details of the rest of the proofs.

Let $\set{Q_k}_{k=1}^\infty$ be the covering of $\Real^d$ by critical balls whose existence is guaranteed by Lemma \ref{bolas}. Consider the auxiliary operator given by
$$f\longmapsto Sf(x)=\sum_{k\geq1}\chi_{Q_k}(x) t\partial_t P_t(\chi_{2Q_k} f)(x),\qquad x\in\Real^d,~t>0.$$
Then $S$ is a bounded operator from $L^p_X(\Real^d)$ into $L^p_{L^q_X((0,\infty),\frac{dt}{t})}(\Real^d)$. Indeed, by using Minkowski's inequality, the finite overlapping of the balls $Q_k$, the boundedness in $L^p$ of $g^{\Delta,q}$ and once more the finite overlapping of $2Q_k$ we get
\begin{align*}
    \norm{Sf}_{L^p_{L^q_X((0,\infty),\frac{dt}{t})}(\Real^d)} &\leq \Big(\int_{\Real^d}\Big|\sum_{k\geq1}\chi_{Q_k}(x)\norm{t\partial_tP_t (\chi_{2Q_k}f)(x)}_{L^q_X((0,\infty),\frac{dt}{t})}\Big|^pdx\Big)^{1/p} \\
    &\leq \sum_{k\geq1}\norm{\chi_{Q_k}g^{\Delta,q}(\chi_{2Q_k}f)}_{L^p_X(\Real^d)}\leq C\sum_{k\geq1}\norm{g^{\Delta,q}(\chi_{2Q_k}f)}_{L^p_X(\Real^d)} \\
    &\leq C\sum_{k\geq1}\norm{\chi_{2Q_k}f}_{L^p_X(\Real^d)}\leq C\norm{f}_{L^p_X(\Real^d)}.
\end{align*}
Recall that for a compactly supported function $f$ in $L^\infty_X(\Real^d)$ we have, as in \eqref{eq:local control},
$$g^{\Delta,q}_{\mathrm{loc}}f(x)\leq g^{\Delta,q}(\chi_N(x,\cdot)f(\cdot))(x)= \norm{t\partial_tP_t\left(\chi_{N}(x,\cdot)f(\cdot)\right)(x)}_{L^q_X((0,\infty),\frac{dt}{t})}, \qquad\hbox{a.e.}~x\in\Real^d.$$
Our idea is to compare the operators $S$ and $f\longmapsto t\partial_tP_t(\chi_{N}(x,\cdot)f(\cdot))(x)$. In order to do that we need some geometrical considerations. Let $C_1$ be the constant that appears in Lemma \ref{Lem:equiv rho}. Consider the set
$$\widetilde{N}=\Big\{(x,y)\in\Real^d\times\Real^d:\abs{x-y}<\frac{C_1}{1+C_1}~\rho(x)\Big\}.$$
It is an exercise to prove that if $(x,y)\in\widetilde{N}$ then, since the family $\set{Q_k}_{k=1}^\infty$ is a covering of $\Real^d$, there exists a positive integer $k$ such that $(x,y)\in Q_k\times2Q_k$. On the other hand, if $(x,y)\in Q_k\times2Q_k$, then by using Lemma \ref{Lem:equiv rho} we get $\abs{x-y}\leq\abs{x-x_k}+\abs{x_k-y}\leq3C_1^{-1}\rho(x)$. Observe that it follows from the finite overlapping property of the balls $Q_k$ that
$$\norm{t\partial_tP_t\left(\chi_N(x,\cdot)f(\cdot)\right)(x)}_X\sim \Big\|\sum_{k\geq1}\chi_{Q_k}(x)t\partial_tP_t\left(\chi_N(x,\cdot)f(\cdot)\right)(x)\Big\|_X,\qquad x\in\Real^d,~t>0.$$
The geometrical comments just made ensure that the kernel of the difference operator
\begin{equation}\label{dif oper}
f\longmapsto\sum_{k\geq1}\chi_{Q_k}(x)t\partial_tP_t\left(\chi_N(x,\cdot)f(\cdot)\right)(x)-Sf(x),\qquad x\in\Real^d,~t>0,
\end{equation}
is supported in the region $\displaystyle A:=\set{(x,y)\in\Real^d\times\Real^d:\frac{C_1}{1+C_1}~\rho(x)\leq\abs{x-y}\leq3C_1^{-1}\rho(x)}$. Consequently, as
\begin{equation}\label{est Pt}
\norm{t\partial_tP_t(x-y)}_{L^q((0,\infty),\frac{dt}{t})}=\frac{C}{\abs{x-y}^d},\qquad x,y\in\Real^d,
\end{equation}
we have
$$\Big\|\sum_{k\geq1}\chi_{Q_k}(x)t\partial_tP_t\left(\chi_N(x,\cdot)f(\cdot)\right)(x)-Sf(x)\Big\|_{L^q((0,\infty),\frac{dt}{t})}\leq C\int_{\Real^d}\frac{\chi_{A}(x,y)}{\abs{x-y}^d}\norm{f(y)}_X~dy.$$
Note that
$$\int_{\Real^d}\frac{\chi_{A}(x,y)}{\abs{x-y}^d}~dy=\int_{\frac{C_1}{1+C_1}\rho(x)\leq\abs{x-y}\leq3C_1^{-1}\rho(x)}\frac{1}{\abs{x-y}^d}~dy= C\log\frac{3(1+C_1)}{C_1^2},\qquad x\in\Real^d,$$
and, by Lemma 2.1,
$$\int_{\Real^d}\frac{\chi_A(x,y)}{\abs{x-y}^d}~dx\leq \int_{\alpha_1\rho(y)\leq\abs{x-y}\leq\alpha_2\rho(y)}\frac{1}{\abs{x-y}^d}~dx=C\log\frac{\alpha_2}{\alpha_1},\qquad y\in\Real^d,$$
for some constants $\alpha_1$ and $\alpha_2$ independent of $y$. Therefore the operator $f\displaystyle\longmapsto \int_{\Real^d}\frac{\chi_{A}(x,y)}{\abs{x-y}^d}\norm{f(y)}_Xdy$ is bounded from $L^p_X$ into $L^p$ for every $p$, $1\leq p<\infty$. Hence we get the conclusion.
\end{proof}

\begin{rem}\label{metodo1}
Consider two Banach spaces $X_1$ and $X_2$. Let $T$ be a linear operator that maps $C^\infty_c(\Real^d;X_1)$ into $X_2$-valued strongly measurable functions. Suppose $T$ has an associated kernel which satisfies the standard Calder\'on-Zygmund estimates. Define the ``$\rho$-localized'' operator
$$T_{\mathrm{loc}}f(x)=T\left(\chi_N(x,\cdot)f(\cdot)\right)(x),\qquad x\in\Real^d,$$
where $N$ is the region determined by $\abs{x-y}\leq\rho(x)$ as above. Then:
\begin{itemize}
\item Assume $T$ has a bounded extension from $L^p_{X_1}(\Real^d)$ into $L^p_{X_2}(\Real^d)$ for some $p$, $1<p<\infty$. Then $T_{\mathrm{loc}}$ has a bounded extension from $L^p_{X_1}(\Real^d)$ into $L^p_{X_2}(\Real^d)$.
\item Assume $T$ has a bounded extension from $L^1_{X_1}(\Real^d)$ into weak-$L^1_{X_2}(\Real^d)$. Then $T_{\mathrm{loc}}$ has a bounded extension from $L^1_{X_1}(\Real^d)$ into weak-$L^1_{X_2}(\Real^d)$.
\item Assume that for every function $f\in L^1_{X_1}(\Real^d)$ we have $\norm{Tf(x)}_{X_2}<\infty$ for almost all $x\in\Real^d$. Then the same is true for $T_{\mathrm{loc}}$.
\end{itemize}
The reader can check the validity of this Remark just by exchanging $X$ by $X_1$, $L^q_X((0,\infty),\frac{dt}{t})$ by $X_2$ and $f\longmapsto t\partial_tP_tf(x)$ by $f\longmapsto Tf(x)$ along the lines of the proof of Theorem \ref{herencia} above.
\end{rem}

The next Theorem permits us to pass, for $\rho$-localized operators related to $-\Delta$, from $L^p$-boundedness to $BMO_{\L}$ and $H^1_{\L}-L^1$ boundedness.

\begin{thm}\label{trasplante}
Let $X$ be a Banach space such that the operator
$$f\longmapsto Tf(x)=t\partial_tP_t\left(\chi_N(x,\cdot)f(\cdot)\right)(x),\qquad x\in\Real^d,~t>0,$$
is bounded from $L^p_X(\Real^d)$ into $L^p_{L^q_X((0,\infty),\frac{dt}{t})}(\Real^d)$ for some $p$, $1<p<\infty$. Then $T$ maps $BMO_{\L,X}$ into $BMO_{\L,L^q_X((0,\infty),\frac{dt}{t})}$ and $H^1_{\L,X}$ into $L^1_{L^q_X((0,\infty),\frac{dt}{t})}(\Real^d)$.
\end{thm}

\begin{proof}
\noindent\textbf{Boundedness from $BMO_{\L,X}$ into $BMO_{\L,L^q_X((0,\infty),\frac{dt}{t})}$.} We first analyze the behavior over ``small'' balls. Consider a ball $B=B(x_0,r_0)$, such that $5r_0<C_1\rho(x_0)$, where $C_1<1$ is the constant that appears in Lemma \ref{Lem:equiv rho}. Given a function $f$ we decompose it as
$$f=(f-f_B)\chi_{4B}+(f- f_B)\chi_{(4B)^c}+f_B=:f_1+f_2+f_3.$$

Before entering into the concrete proof, we need some small preparation. For $x,z\in B$,
$$Tf(x)-Tf(z)=Tf_1(x)-Tf_1(z)+Tf_2(x)-Tf_2(z)+Tf_3(x)-Tf_3(z).$$
We begin by observing that
\begin{align*}
    Tf_2(x)-Tf_2(z)+Tf_3(x)-Tf_3(z) &= \int_{\Real^d}\left(t\partial_tP_t(x-y)-t\partial_tP_t(z-y)\right)\chi_{\abs{x-y}\leq\rho(x)}(y)f_2(y)~dy \\
    &\quad +\int_{\Real^d}t\partial_tP_t(z-y)\left(\chi_{\abs{x-y}\leq\rho(x)}(y)-\chi_{\abs{z-y}\leq\rho(x)}(y)\right)f_2(y)~dy \\
    &\quad +\int_{\Real^d}t\partial_tP_t(z-y)\left(\chi_{\abs{z-y}\leq\rho(x)}(y)-\chi_{\abs{z-y}\leq\rho(z)}(y)\right)f_2(y)~dy \\
    &\quad +f_B\int_{\Real^d}\left(t\partial_tP_t(x-y)-t\partial_tP_t(z-y)\right)\chi_{\abs{x-y}\leq\rho(x)}(y)~dy \\
    &\quad +f_B\int_{\Real^d}t\partial_tP_t(z-y)\left(\chi_{\abs{x-y}\leq\rho(x)}(y)-\chi_{\abs{z-y}\leq\rho(x)}(y)\right)~dy \\
    &\quad +f_B\int_{\Real^d}t\partial_tP_t(z-y)\left(\chi_{\abs{z-y}\leq\rho(x)}(y)-\chi_{\abs{z-y}\leq\rho(z)}(y)\right)~dy. \\
\end{align*}
Using Lemma \ref{Lem:equiv rho}, $\chi_{(4B)^c}(y)\left(\chi_{\abs{z-y}\leq\rho(x)}(y)-\chi_{\abs{z-y}\leq\rho(z)}(y)\right)= \chi_{\abs{z-y}\leq\rho(x)}(y)-\chi_{\abs{z-y}\leq\rho(z)}(y)$. Therefore,
\begin{multline*}
\int_{\Real^d}t\partial_tP_t(z-y)\left(\chi_{\abs{z-y}\leq\rho(x)}(y)-\chi_{\abs{z-y}\leq\rho(z)}(y)\right)f_2(y)~dy \\
=\int_{\Real^d}t\partial_tP_t(z-y)\left(\chi_{\abs{z-y}\leq\rho(x)}(y)-\chi_{\abs{z-y}\leq\rho(z)}(y)\right)(f(y)-f_B)~dy.
\end{multline*}
As a consequence,
\begin{align*}
    Tf_2(x)-Tf_2(z)+Tf_3(x)-Tf_3(z) &= \int_{\Real^d}\left(t\partial_tP_t(x-y)-t\partial_tP_t(z-y)\right)\chi_{\abs{x-y}\leq\rho(x)}(y)f_2(y)~dy \\
    &\quad +\int_{\Real^d}t\partial_tP_t(z-y)\left(\chi_{\abs{x-y}\leq\rho(x)}(y)-\chi_{\abs{z-y}\leq\rho(x)}(y)\right)f_2(y)~dy \\
    &\quad +\int_{\Real^d}t\partial_tP_t(z-y)\left(\chi_{\abs{z-y}\leq\rho(x)}(y)-\chi_{\abs{z-y}\leq\rho(z)}(y)\right)f(y)~dy \\
    &\quad +f_B\int_{\Real^d}\left(t\partial_tP_t(x-y)-t\partial_tP_t(z-y)\right)\chi_{\abs{x-y}\leq\rho(x)}(y)~dy \\
    &\quad +f_B\int_{\Real^d}t\partial_tP_t(z-y)\left(\chi_{\abs{x-y}\leq\rho(x)}(y)-\chi_{\abs{z-y}\leq\rho(x)}(y)\right)~dy \\
    &=: A_1(x,z)+ A_2(x,z)+A_3(x,z)+A_4(x,z)+A_5(x,z).
\end{align*}

After these remarks, we can start the actual proof of the boundedness in $BMO$. We have
\begin{multline*}
\frac{1}{\abs{B}}\int_B\norm{Tf(x)-(Tf)_B}_{L^q_X((0,\infty),\frac{dt}{t})}dx \\
\leq\frac{2}{\abs{B}}\int_B\norm{Tf_1(x)}_{L^q_X((0,\infty),\frac{dt}{t})}dx+\sum_{i=1}^5\frac{1}{\abs{B}^2} \int_B\int_B\norm{A_i(x,z)}_{L^q_X((0,\infty),\frac{dt}{t})}dx~dz.
\end{multline*}
By hypothesis $T$ is bounded from $L^p_X(\Real^d)$ into $L^p_{L^q_X((0,\infty),\frac{dt}{t})}(\Real^d)$, so
\begin{align*}
    \frac{1}{\abs{B}}\int_B\norm{Tf_1(x)}_{L^q_X((0,\infty),\frac{dt}{t})}dx &\leq C\left(\frac{1}{\abs{B}}\int_B\norm{Tf_1(x)}^p_{L^q_X((0,\infty),\frac{dt}{t})}dx\right)^{1/p} \\
    &\leq C\left(\frac{1}{\abs{B}}\int_{\Real^d}\norm{f_1(x)}^p_X~dx\right)^{1/p} \\
    &= C\left(\frac{1}{\abs{B}}\int_{4B}\norm{f(x)-f_B}^p_X~dx\right)^{1/p} \\
    &\leq C\norm{f}_{BMO_X}\leq C\norm{f}_{BMO_{\L,X}},
\end{align*}
where in the penultimate inequality we applied an argument as in \eqref{argument} below. Let us now estimate all the $A_i(x,z)$, $i=1,\ldots,5$, for $x,z\in B=B(x_0,r_0)$. By the Mean Value Theorem and \eqref{est Pt},
\begin{align*}
    \norm{A_1(x,z)}_{L^q_X((0,\infty),\frac{dt}{t})} &\leq C\int_{\Real^d}\frac{\abs{x-z}}{\abs{x-y}^{d+1}}\norm{f_2(y)}_Xdy \\
     &\leq Cr_0\int_{\abs{x_0-y}>4r_0}\frac{1}{\abs{x_0-y}^{d+1}}\norm{f(y)-f_B}_Xdy \\
     &= Cr_0\sum_{j=2}^\infty\int_{2^jr_0<\abs{x_0-y}\leq2^{j+1}r_0}\frac{1}{\abs{x_0-y}^{d+1}}\norm{f(y)-f_B}_Xdy \\
     &\leq C\sum_{j=0}^\infty\frac{1}{2^j}\frac{1}{(2^{j+1}r_0)^d}\int_{\abs{x_0-y}\leq2^{j+1}r_0}\norm{f(y)-f_B}_Xdy \\
     &\leq C\norm{f}_{BMO_X}\leq C\norm{f}_{BMO_{X,\L}}.
\end{align*}
Again by \eqref{est Pt},
$$\norm{A_2(x,z)}_{L^q_X((0,\infty),\frac{dt}{t})}\leq C\int_{\Real^d}\frac{1}{\abs{z-y}^d} \abs{\chi_{\abs{x-y}\leq\rho(x)}(y)-\chi_{\abs{z-y}\leq\rho(x)}(y)}\norm{f_2(y)}_Xdy.$$
Observe that $A_2$ will be non zero in the following cases:
\begin{enumerate}[(i)]
    \item $\abs{x-y}\leq\rho(x)$ and $\abs{z-y}>\rho(x)$,
    \item $\abs{x-y}>\rho(x)$ and $\abs{z-y}\leq\rho(x)$.
\end{enumerate}
In the first case $\rho(x)<\abs{z-y}\leq\abs{z-x}+\abs{x-y}<2r_0+\abs{x-y}$ and then $\rho(x)-2r_0<\abs{x-y}\leq\rho(x)$. While in $\mathrm{(ii)}$ we have $\rho(x)<\abs{x-y}\leq\abs{x-z}+\abs{z-y}<2r_0+\rho(x)$. On the other hand $\abs{x-y}\sim\abs{z-y}$. Let $j_0$ and $j_1$ be nonnegative integers such that $2^{j_0}r_0\leq\rho(x)/2<2^{j_0+1}r_0$ and $2^{j_1}r_0\leq2\rho(x)<2^{j_1+1}r_0$. Observe that, since $5r_0<\rho(x)$ for all $x\in B(x_0,r_0)$, we have $j_0\geq1$. The Mean Value Theorem gives $(\rho(x)-2r_0)^d-(\rho(x)+2r_0)^d\leq C\rho(x)^{d-1}r_0$, hence applying H\"older's inequality with some $r\in(1,\infty)$ we get
\begin{align*}
    \lefteqn{\norm{A_2(x,z)}_{L^q_X((0,\infty),\frac{dt}{t})} \leq C\int_{\rho(x)-2r_0<\abs{x-y}<\rho(x)+2r_0}\frac{1}{\abs{x-y}^d}\norm{f_2(y)}_Xdy} \\
    &\leq C\left(\int_{\rho(x)-2r_0<\abs{x-y}<\rho(x)+2r_0}\frac{1}{\abs{x-y}^{dr}} \norm{f(y)-f_B}_X^rdy\right)^{1/r}\rho(x)^{(d-1)/r'}r_0^{1/r'} \\
    &\leq C\left(\int_{\rho(x)/2<\abs{x-y}<2\rho(x)}\frac{1}{\abs{x-y}^{dr}} \norm{f(y)-f_B}_X^rdy\right)^{1/r}\rho(x)^{(d-1)/r'}r_0^{1/r'} \\
    &\leq C\left(\sum_{j=j_0}^{j_1}\frac{1}{(2^jr_0)^{(d-1)(r-1)}}\frac{1}{2^{j(r-1)}}\frac{1}{(2^jr_0)^d} \int_{\abs{x_0-y}<2^{j+2}r_0}\norm{f(y)-f_B}_X^rdy\right)^{1/r}\rho(x)^{(d-1)/r'} \\
    &\leq C\left(\frac{1}{(2^{j_0}r_0)^{(d-1)(r-1)}}\sum_{j=j_0}^{j_0+2}\frac{1}{2^{j(r-1)}} \frac{1}{(2^jr_0)^d}\int_{\abs{x_0-y}<2^{j+2}r_0}\norm{f(y)-f_B}_X^rdy\right)^{1/r}\rho(x)^{(d-1)/r'} \\
    &\leq C\left(\sum_{j=0}^\infty\frac{1}{2^{j(r-1)}}\frac1{(2^jr_0)^d} \int_{\abs{x_0-y}<2^{j+2}r_0}\norm{f(y)-f_B}_X^rdy\right)^{1/r}\leq C\norm{f}_{BMO_X}.
\end{align*}
Observe that in the penultimate inequality above we pass to the infinite series since $j_0$ depends on $\rho(x)$ and we want an estimate independent of it. For the last inequality above we first note that, by the triangle inequality and Minkowski's integral inequality,
\begin{align}
\nonumber    & \left(\frac1{(2^j r_0)^d}\int_{\abs{x_0-y}<2^{j+2}r_0}\norm{f(y)-f_B}_X^rdy\right)^{1/r} \\
\label{argument}     & \leq\left(\frac{1}{(2^jr_0)^d}\int_{\abs{x_0-y}<2^{j+2}r_0}\left(\norm{f(y)-f_{2^{j+2}B}}_X+\sum_{k=0}^{j+1} \norm{f_{2^{k+1}B}-f_{2^kB}}_X\right)^rdy\right)^{1/r} \\
\nonumber     & \leq\left(\frac{4^d}{(2^{j+2}r_0)^d} \int_{\abs{x_0-y}<2^{j+2}r_0} \norm{f(y)-f_{2^{j+2}B}}^r_Xdy\right)^{1/r}+C\sum_{k=0}^{j+1}\norm{f_{2^{k+1}B}-f_{2^kB}}_X \\
\nonumber     & \leq C(j+3)\norm{f}_{BMO_X}.
\end{align}
Hence,
$$\sum_{j=0}^\infty\frac{1}{2^{j(r-1)}}\frac1{(2^j r_0)^d} \int_{\abs{x_0-y}<2^{j+2}r_0}\norm{f(y)-f_B}_X^rdy\leq  C\sum_{j=0}^\infty\frac{(j+3)^r}{2^{j(r-1)}}\norm{f}^r_{BMO_X}=C\norm{f}^r_{BMO_X}.$$
Since $x,z\in B$,
\begin{multline*}
    \norm{A_3(x,z)}_{L^q_X((0,\infty),\frac{dt}{t})} \leq \int_{C_1\rho(x_0)<\abs{z-y}<C_1^{-1}\rho(x_0)}\frac{1}{\abs{z-y}^d}\norm{f(y)}_Xdy \\
    \leq \frac{C}{\rho(x_0)^d}\int_{\abs{z-y}<C_1^{-1}\rho(x_0)}\norm{f(y)}_Xdy
    \leq \frac{C}{\rho(z)^d}\int_{\abs{z-y}<C\rho(z)}\norm{f(y)}_Xdy\leq C\norm{f}_{BMO_{\L,X}}.
\end{multline*}
By dominated convergence,
$$\int_{\Real^d}\partial_tP_t(x-y)~dy=\partial_t\int_{\Real^d}P_t(x-y)~dy=\partial_t1=0,\qquad x\in\Real^d.$$
Therefore,
\begin{align*}
    \norm{A_4(x,z)}_{L^q_X((0,\infty),\frac{dt}{t})} &= \norm{f_B\int_{\Real^d}\left(t\partial_tP_t(x-y) -t\partial_tP_t(z-y)\right)\chi_{\abs{x-y}\leq\rho(x)}(y)~dy}_{L^q_X((0,\infty),\frac{dt}{t})} \\
    &= \Big\|f_B\int_{\Real^d}\left(t\partial_tP_t(x-y)-t\partial_tP_t(z-y)\right)\chi_{\abs{x-y}\leq\rho(x)}(y)~dy \\
    &\qquad -f_B\int_{\Real^d}\left(t\partial_tP_t(x-y)-t\partial_tP_t(z-y)\right)~dy\Big\|_{L^q_X((0,\infty),\frac{dt}{t})} \\
    &= \norm{f_B\int_{\Real^d}\left(t\partial_tP_t(x-y) -t\partial_tP_t(z-y)\right)\chi_{\abs{x-y}>\rho(x)}(y)~dy}_{L^q_X((0,\infty),\frac{dt}{t})} \\
    &\leq C\norm{f_B}_X\int_{\Real^d}\Big\|t\partial_tP_t(x-y)-t\partial_tP_t(z-y)\Big\|_{L^q((0,\infty)\frac{dt}{t})}\chi_{\abs{x-y}>\rho(x)}(y)~dy \\
    &\leq C\norm{f_B}_X\int_{\Real^d}\frac{\abs{x-z}}{\abs{x-y}^{d+1}}\chi_{\abs{x-y}>\rho(x)}(y)~dy \\
    &\leq C\norm{f_B}_X\frac{r_0}{\rho(x)}\leq C\norm{f_B}_X\frac{r_0}{\rho(x_0)}.
\end{align*}
As $\norm{f_B}_X\leq C\left(1+\log\frac{\rho(x_0)}{r_0}\right)\norm{f}_{BMO_{\L,X}}$ (see \cite[Lemma~2]{DGMTZ}) we get the appropriate bound for $A_4$. Finally, by using the arguments in $A_2$,
\begin{align*}
    \norm{A_5(x,z)}_{L^q_X((0,\infty),\frac{dt}{t})} &\leq C\norm{f_B}_X\int_{\rho(x)-2r_0<\abs{x-y}<\rho(x)+2r_0} \frac{1}{\abs{x-y}^d}~dy \\
    &\leq C\norm{f}_{BMO_{\L,X}}\left(1+\log\frac{\rho(x_0)}{r_0}\right)\log\left(\frac{\rho(x)+2r_0}{\rho(x)-2r_0}\right)
\end{align*}
Since $\frac{r_0}{\rho(x)}<1/5$, we have $\log\left(\frac{\rho(x)+2r_0}{\rho(x)-2r_0}\right)\sim\frac{r_0}{\rho(x)}\sim\frac{r_0}{\rho(x_0)}$, that gives the desired bound for $A_5$.

Let us now analyze the behavior over ``big'' balls. Let $B_1=B(x_0,k\rho(x_0))$, with $k\geq\frac{C_1}{5}$. Given a function $f$ we decompose it as $f=f_1+f_2$, where $f_1=f\chi_{2B_1}$. By H\"older's inequality and the hypothesis,
$$\frac{1}{\abs{B_1}}\int_{B_1}\norm{t\partial_tP_t\left(\chi_N(x,\cdot)f_1(\cdot)\right)(x)}_{L^q_X((0,\infty), \frac{dt}{t})}dx\leq C\left(\frac{1}{\abs{B_1}}\int_{2B_1}\norm{f(x)}^p_Xdx\right)^{1/p}\leq C\norm{f}_{BMO_{\L,X}}.$$
On the other hand, by using Lemma \ref{Lem:equiv rho} and \eqref{est Pt},
\begin{align*}
    \norm{t\partial_tP_t\left(\chi_N(x,\cdot)f_2(\cdot)\right)(x)}_{L^q_X((0,\infty),\frac{dt}{t})} &\leq \int_{\Real^d}\frac{1}{\abs{x-y}^d}\chi_{2k\rho(x_0)\leq\abs{x-y}\leq\rho(x)}(y)\norm{f(y)}_Xdy \\
    &\leq \frac{C}{\rho(x)^d}\int_{\abs{x-y}\leq\rho(x)}\norm{f(y)}_Xdy\leq C\norm{f}_{BMO_{\L,X}},\quad x\in B_1.
\end{align*}
This finishes the proof of the $BMO$ boundedness.

\noindent\textbf{Boundedness from $H^1_{\L,X}$ into $L^1_{L^q_X((0,\infty),\frac{dt}{t})}(\Real^d)$.} We begin with the analysis over atoms supported on ``small'' balls. Let $a$ be an atom with support contained in a ball $\tilde{B}=B(y_0,r_0)$, with $r_0<\rho(y_0)$. Then
$$\int_{\Real^d}\norm{Ta(x)}_{L^q_X((0,\infty),\frac{dt}{t})}dx=\int_{4\tilde{B}}\norm{Ta(x)}_ {L^q_X((0,\infty),\frac{dt}{t})}dx+\int_{(4\tilde{B})^c}\norm{Ta(x)}_{L^q_X((0,\infty),\frac{dt}{t})}dx=:A_1+A_2.$$
Since $T$ is bounded in $L^p$, by \eqref{prop atom} we have
$$A_1\leq C\left(\int_{4\tilde{B}}\norm{Ta(x)}^p_{L^q_X((0,\infty),\frac{dt}{t})}dx\right)^{1/p}|\tilde{B}|^{1/p'}\leq C\left(\int_{\tilde{B}}\norm{a(x)}^p_Xdx\right)^{1/p}|\tilde{B}|^{1/p'}\leq C.$$
Applying the fact that the atom $a$ has mean zero \eqref{prop atom media} we get
\begin{align*}
    A_2 &= \int_{(4\tilde{B})^c}\norm{\int_{\Real^d}\left(t\partial_t P_t(x-y)\chi_{\abs{x-y}\leq\rho(x)}(y)-t\partial_tP_t(x-y_0) \chi_{\abs{x-y_0}\leq\rho(x)}(x)\right)a(y)~dy}_{L^q_X((0,\infty),\frac{dt}{t})}dx \\
    &\leq \int_{(4\tilde{B})^c}\norm{\int_{\Real^d}\left(t\partial_tP_t(x-y)-t\partial_tP_t(x-y_0)\right) \chi_{\abs{x-y}\leq\rho(x)}(y)a(y)~dy}_{L^q_X((0,\infty),\frac{dt}{t})}dx \\
    &\quad +\int_{(4\tilde{B})^c}\norm{\int_{\Real^d}t\partial_tP_t(x-y_0) \left(\chi_{\abs{x-y}\leq\rho(x)}(y)-\chi_{\abs{x-y_0}\leq\rho(x)}(x)\right)a(y)~dy}_{L^q_X((0,\infty),\frac{dt}{t})}dx \\
    &\leq C\int_{(4\tilde{B})^c}\int_{\Real^d}\frac{\abs{y-y_0}}{\abs{x-y_0}^{d+1}}\norm{a(y)}_Xdy~dx \\
    &\quad +C\int_{(4\tilde{B})^c}\int_{\Real^d}\frac{1}{\abs{x-y_0}^d}\abs{\chi_{\abs{x-y}\leq\rho(x)}(y)- \chi_{\abs{x-y_0}\leq\rho(x)}(x)}\norm{a(y)}_Xdy~dx \\
    &=: C\left(A_{21}+A_{22}\right).
\end{align*}
Fubini's Theorem and \eqref{prop atom} give
\begin{align*}
    A_{21} &= \int_{\Real^d}\abs{y-y_0}\norm{a(y)}_X\left[\int_{\Real^d}\chi_{\abs{x-y_0}\geq4r_0}(x) \frac{1}{\abs{x-y_0}^{d+1}}~dx\right]dy \\
    &= \frac{C}{r_0}\int_{\abs{y-y_0}<r_0}\abs{y-y_0}\norm{a(y)}_Xdy\leq C.
\end{align*}
A geometric reasoning parallel to the one developed above for the $BMO$ case gives that in order to $A_{22}\neq0$ we must have $3r_0<\rho(x)$, $\rho(x)-r_0<\abs{x-y}<\rho(x)+r_0$ and, in addition, $\abs{x-y_0}\sim\rho(x)\sim\rho(y_0)$. Therefore, since the atom $a$ is supported in $\tilde{B}=B(y_0,r_0)$ and is controlled in $L^\infty$ norm by $Cr_0^{-d}$,
\begin{align*}
    A_{22} &\leq \frac{C}{\rho(y_0)^d}\int_{\abs{x-y_0}\leq C\rho(y_0)}\int_{\rho(x)-r_0<\abs{x-y}<\rho(x)+r_0}\norm{a(y)}_Xdy~dx \\
     &\leq \frac{C}{\rho(y_0)^d}\int_{\abs{x-y_0}\leq C\rho(y_0)}\int_{\abs{y-y_0}<r_0}\norm{a(y)}_X~dy~dx\leq\frac{C}{\rho(y_0)^d}\int_{\abs{x-y_0}\leq C\rho(y_0)}dx\leq C.
\end{align*}

We continue with the analysis over atoms supported on ``big'' balls. Let $a$ be an atom supported in a ball $\bar{B}(y_0,\gamma\rho(y_0))$, with $\gamma>1$. We begin by proceeding as in the previous case for $A_1$. For $A_2$, since we do not have the cancelation property \eqref{prop atom media}, we estimate its size as follows:
$$A_2=\int_{(4\bar{B})^c}\norm{Ta(x)}_{L^q_X((0,\infty)\frac{dt}{t})}dx\leq C\int_{(4\bar{B})^c}\int_{\bar{B}}\frac{1}{\abs{x-y}^d}~\norm{a(y)}_X\chi_{\abs{x-y}\leq\rho(x)}(y)~dy~dx.$$
The domain of integration above is contained in the set defined by the conditions $\abs{x-y_0}\geq4\gamma\rho(y_0)$, $\abs{y-y_0}<\gamma\rho(y_0)$ and $\abs{x-y}\leq\rho(x)$. These conditions imply that $4\gamma\rho(y_0)\leq\abs{x-y_0}\leq\abs{x-y}+\abs{y-y_0}<\abs{x-y}+\gamma\rho(y_0)$, hence $3\gamma\rho(y_0)\leq\abs{x-y}$. Note that, by Lemma \ref{Lem:equiv rho}, $\rho(x)\leq C\rho(y)\leq\bar{C}\gamma\rho(y_0)$. Therefore $3\gamma\rho(y_0)\leq\abs{x-y}\leq\bar{C}\gamma\rho(y_0)$ and we get
$$A_2\leq C\int_{\bar{B}}\norm{a(y)}_X\int_{3\gamma\rho(y_0)\leq\abs{x-y}\leq\bar{C}\gamma\rho(y_0)}\frac{1}{\abs{x-y}^d}~dx~dy\leq C.$$
\end{proof}

\begin{rem}\label{metodo2}
Consider two Banach spaces $X_1$ and $X_2$. Let $T$ be a linear operator that maps $L^p_{X_1}(\Real^d)$ into $L^p_{X_2}(\Real^d)$ for some $p$, $1<p<\infty$, such that $T1$ can be defined and $T1=0$. Assume $T$ has an associated kernel which satisfies the standard estimates of Calder\'on-Zygmund operators. Define the operator
$$T_{\mathrm{loc}}f(x)=T\left(\chi_N(x,\cdot)f(\cdot)\right)(x),\qquad x\in\Real^d.$$
Then:
\begin{itemize}
    \item $T_{\mathrm{loc}}$ is bounded from $BMO_{\L,X_1}$ into $BMO_{\L,X_2}$, and
    \item $T_{\mathrm{loc}}$ is bounded from $H^1_{\L,X_1}$ into $L^1_{X_2}(\Real^d)$.
\end{itemize}
Parallel to Remark \ref{metodo1}, the reader can check the validity of these claims just by exchanging, along the lines of the proof of Theorem \ref{trasplante}, $X$ by $X_1$, $L^q_X((0,\infty),\frac{dt}{t})$ by $X_2$ and $f\longmapsto t\partial_tP_tf(x)$ by $f\longmapsto Tf(x)$.
\end{rem}

\section{Proof of Theorem A}\label{Section:Proof}

Given a Banach space $X$, define the modulus of convexity by
$$\delta_X(\varepsilon) = \inf\set{1-\norm{\frac{x+y}{2}}:x,y\in X,\norm{x}=\norm{y}=1,\norm{x-y}=\varepsilon},\quad0<\varepsilon<2.$$
The Banach space $X$ is called $q$-uniformly convex, $2\leq q<\infty$, if $\delta_X(\varepsilon)\geq c\varepsilon^q$ for some positive constant $c$. By Pisier's Renorming Theorem \cite{Pisier}, $X$ is $q$-uniformly convex if and only if $X$ is of martingale cotype $q$. For martingale cotype the following Theorem holds, see \cite{Xu} and \cite{Martinez-Torrea-Xu}.

\begin{thm}\label{maite}
Let $X$ be a Banach space and $2\leq q<\infty$.  The following statements are equivalent.
\begin{enumerate}[(1)]
    \item $X$ is of martingale cotype $q$.
    \item The operator $g^{\Delta,q}$ maps $BMO_{c,X}$ into $BMO$.
    \item The operator $g^{\Delta,q}$ maps $L^p_X(\Real^d)$ into $L^p(\Real^d)$, for any $p$ in the range $1<p<\infty$.
    \item The operator $g^{\Delta,q}$ maps $L^1_X(\Real^d)$ into weak-$L^1(\Real^d)$.
    \item The operator $g^{\Delta,q}$ maps $H^1_{X}$ into $L^1(\Real^d)$.
    \item For every $f\in L^1_X(\Real^d)$, $g^{\Delta,q}f(x)<\infty$ for almost every $x\in\Real^d$.
\end{enumerate}
\end{thm}

The space $H^1_X$ denotes the atomic Hardy space in $\Real^d$. By $BMO_{c,X}$ we mean the set of functions that belong to the classical $BMO$ with values in $X$ and have compact support.

\begin{proof}[Proof of Theorem A]
Observe that hypothesis $(i)$ is equivalent to one of the statements in Theorem \ref{maite}.

$(i)\implies(ii)$.  We can apply Theorems \ref{herencia} and \ref{trasplante} to get that the operator $f\longmapsto t\partial_tP_t\left(\chi_N(x,\cdot)f(\cdot)\right)$ maps $BMO_{\L,X}$ into $BMO_{\L,L^q_X((0,\infty),\frac{dt}{t})}$. By using Lemma \ref{difloc} we obtain the boundedness from $BMO_{\L,X}$ into $BMO_{\L}$ of the operator $g^{\L,q}_{\mathrm{loc}}$. Finally, by Lemma \ref{Lem:ltodop} $\mathrm{(b)}$ we arrive to $(ii)$.

$(i)\implies(iii)$. By Theorem  \ref{herencia} and Lemma \ref{difloc} $\mathrm{(a)}$ the local operator $g^{\L,q}_{\mathrm{loc}}$ is bounded in $L^p$. Boundedness of the global part follows from Lemma \ref{Lem:ltodop} $\mathrm{(a)}$.

$(i)\implies(iv)$. Theorem  \ref{herencia} and Lemma \ref{difloc} $\mathrm{(a)}$, together with Lemma \ref{Lem:ltodop} $\mathrm{(a)}$, give the conclusion.

$(i)\implies(v)$. By using Theorems \ref{herencia}, \ref{trasplante} and \ref{difloc} $\mathrm{(c)}$ we see that $g^{\L,q}_{\mathrm{loc}}$ maps $H^1_{\L,X}$ into $L^1(\Real^d)$. Then Lemma \ref{Lem:ltodop} $\mathrm{(c)}$ gives the result.

$(i)\implies(vi)$. Apply Theorem \ref{herencia} and Lemmas \ref{difloc} $\mathrm{(a)}$ and \ref{Lem:ltodop} $\mathrm{(a)}$.

$(ii)\implies(i)$. Theorem \ref{maite} tells us that it is enough to prove the boundedness of $g^{\Delta,q}$ from $BMO_{c,X}$ into $BMO$. From the hypothesis, Lemma \ref{Lem:ltodop} (b) and \eqref{eq:dif con corte} we can deduce that the operator $f(x)\longmapsto g^{\L,q}\left(\chi_N(x,\cdot)f(\cdot)\right)(x)$, $x\in\Real^d$, is bounded from $BMO_{\L,X}$ into $BMO_{\L}$. On the other hand, the proof of Lemma \ref{difloc} shows that the difference operator $f(x)\longmapsto g^{\L,q}\left(\chi_N(x,\cdot)f(\cdot)\right)(x)-g^{\Delta,q}\left(\chi_N(x,\cdot)f(\cdot)\right)(x)$ is bounded from $BMO_{\L,X}$ into $L^\infty$. Thus
$$\hbox{the operator}~f(x)\longmapsto g^{\Delta,q}\left(\chi_N(x,\cdot)f(\cdot)\right)(x)~\hbox{is bounded from}~BMO_{\L,X}~\hbox{into}~BMO_{\L}\subset BMO.$$
Let $f$ be a function in $BMO_{c,X}$. Given a ball $B(x_0,s)$, by Lemma \ref{encoger} there exists $R>0$ depending on $s$ and the support of $f$ such that $\supp f^R\subset B(0,\frac{\rho(0)}{2})$ (see the proof of Lemma \ref{encoger}) and
\begin{align*}
    \frac{1}{\abs{B(x_0,s)}}\int_{B(x_0,s)}g^{\Delta,q}f(x)~dx &= \frac{1}{\abs{B(x_0,s)}}\int_{B(x_0,s)}g^{\Delta,q}\left(\chi_N(\tfrac{x}{R},\cdot)f^R(\cdot)\right)(\tfrac{x}{R})~dx \\
    &= \frac{1}{\abs{B(\tfrac{x_0}{R},\tfrac{s}{R})}}\int_{B(\tfrac{x_0}{R},\tfrac{s}{R})} g^{\Delta,q}(\chi_N(z,\cdot)f^R(\cdot))(z)~dz.
\end{align*}
Since $R$ can be arbitrarily large, we fix it in such a way that $(R\rho(0))^{-d}\norm{f}_{L^1_X(\Real^d)}\leq\norm{f}_{BMO_X}$. Therefore,
\begin{align*}
    \lefteqn{\frac{1}{\abs{B(x_0,s)}}\int_{B(x_0,s)}\abs{g^{\Delta,q}f(x)- \left(g^{\Delta,q}\left(\chi_N(z,\cdot)f^R(\cdot)\right)(z)\right)_{B(\tfrac{x_0}{R},\tfrac{s}{R})}}~dx} \\
    &= \frac{1}{\abs{B(\tfrac{x_0}{R},\tfrac{s}{R})}}\int_{B(\tfrac{x_0}{R},\tfrac{s}{R})}\abs{g^{\Delta,q}\left(\chi_N(x,\cdot)f^R(\cdot)\right)(x) -\left(g^{\Delta,q}\left(\chi_N(z,\cdot)f^R(\cdot)\right)(z)\right)_{B(\tfrac{x_0}{R},\tfrac{s}{R})}}dx \\
    &\leq C\norm{f^R}_{BMO_{\L,X}}\leq C\norm{f}_{BMO_X},
\end{align*}
where for the last inequality above the following argument is applied. Note that to have such an inequality we only have to compare the integral means of $f^R$ with the $BMO_X$-norm of $f$. Let $\alpha\geq1$. If $B(x,\alpha\rho(x))$ does not intersect $B(0,\frac{\rho(0)}{2})$ then $\displaystyle\int_{B(x,\alpha\rho(x))}\norm{f^R(y)}_X~dy=0$ and there is nothing to prove. In case $B(x,\alpha\rho(x))\cap B(0,\frac{\rho(0)}{2})\neq\emptyset$ then, by Lemma \ref{Lem:equiv rho}, $\rho(x)\sim\rho(0)$ and, by the choice of $R$,
\begin{align*}
    \frac{1}{\abs{B(x,\alpha\rho(x))}}\int_{B(x,\alpha\rho(x))}\norm{f^R(y)}_Xdy &\leq \frac{C_n}{(R\alpha\rho(x))^d}\int_{B(0,R\frac{\rho(0)}{2})} \norm{f(z)}_Xdz \\
     & \leq\frac{C}{(R\rho(0))^d}\norm{f}_{L^1_X(\Real^d)}\leq C\norm{f}_{BMO_X};
\end{align*}
here the constant $C$ is independent of $f$.

$(iii)\implies(i)$. Lemmas \ref{Lem:ltodop} $\mathrm{(a)}$ and \ref{difloc} $\mathrm{(a)}$ assure that $g^{\Delta,q}\left(\chi_Nf\right)$ is bounded from $L^p_X(\Real^d)$ into $L^p(\Real^d)$. Let $f\in L^p_X(\Real^d)$ be a function with support contained in a ball $B_M=B(0,M)$, $M>0$. By Lemma \ref{encoger} we can find $R>0$ such that $g^{\Delta,q}f(x)=g^{\Delta,q}\left(\chi_N(\tfrac{x}{R},\cdot)f^R(\cdot)\right)(\tfrac{x}{R})$, for all $\abs{x}<M$. Hence
\begin{multline*}
      \norm{\chi_{B_M}g^{\Delta,q}f}_{L^p(\Real^d)}^p = \int_{\Real^d}\abs{\chi_{B_M}(x)g^{\Delta,q}\left(\chi_N(\tfrac{x}{R},\cdot)f^R(\cdot)\right)(\tfrac{x}{R})}^pdx \\
      \leq R^d\int_{\Real^d}\abs{g^{\Delta,q}\left(\chi_N(\tfrac{x}{R},\cdot)f^R(\cdot)\right)(\tfrac{x}{R})}^pdx\leq CR^d\int_{\Real^d}\norm{f^R(x)}_X^pdx=C\norm{f}_{L^p_X(\Real^d)}^p.
\end{multline*}
As the constant $C$ does not depend on $M$ we can take $M\to\infty$ to get $\norm{g^{\Delta,q}f}_{L^p(\Real^d)}\leq C\norm{f}_{L^p_X(\Real^d)}$.

$(iv)\implies(i)$. We leave this case to the reader.

$(v)\implies(i)$. By Theorem \ref{maite} it is enough to prove the boundedness of $g^{\Delta,q}$ from $H^1_X$ into $L^1(\Real^d)$. Lemmas \ref{Lem:ltodop} $\mathrm{(c)}$ and \ref{difloc} $\mathrm{(c)}$ imply that the localized operator $f\longmapsto\norm{t\partial_tP_t(\chi_N(x,\cdot)f(\cdot))(x)}_{L^q_X((0,\infty),\frac{dt}{t})}$ is bounded from $H^1_{\L,X}$ into $L^1(\Real^d)$. Therefore we only have to prove the boundedness over $H^1$-atoms with cancelation but supported in big balls. Let $a$ be such an atom, namely a function supported in a ball $B(y_0,\gamma\rho(y_0))$ with $\gamma>1$ and $\int_{\Real^d}a(y)~dy=0$. Consider the function $\widetilde{a}^R(x):=R^da^R(x)=R^da(Rx)$, $x\in\Real^d$, $R>0$. The function $\widetilde{a}^R$ is an atom with support contained in the ball $B(\frac{y_0}{R},\frac{\gamma y_0}{R})$. Given $M>0$, Lemma \ref{encoger} allows us to choose a sufficiently large $R$ such that $g^{\Delta,q}a(x)=g^{\Delta,q}\left(\chi_N(\tfrac{x}{R},\cdot)a^R(\cdot)\right)(\tfrac{x}{R})$, for $\abs{x}<M$. Hence
\begin{align*}
    \int_{\abs{x}<M}\abs{g^{\Delta,q}a(x)}~dx &= \int_{\abs{x}<M}\abs{g^{\Delta,q}\left(\chi_N(\tfrac{x}{R},\cdot)a^R(\cdot)\right)(\tfrac{x}{R})}~dx \\
    &= \int_{\abs{z}<\tfrac{M}{R}}\abs{g^{\Delta,q}(\chi_N(z,\cdot)\widetilde{a}^R(\cdot))(z)}~dz \\
    &\leq C\norm{\widetilde{a}^R}_{H^1_{\L,X}}=C\norm{\widetilde{a}^R}_{H^1_{X}}\leq C,
\end{align*}
where $C$ does not depend on $M$. To conclude take $M\to\infty$.

$(vi)\implies(i)$.  We will prove that $g^{\Delta,q}f(x)<\infty$ for almost every $x\in\Real^d$, see Theorem \ref{maite}. By Lemma \ref{Lem:ltodop} $\mathrm{(a)}$ we have that $g^{\L,q}_{\mathrm{loc}}f(x)<\infty$ for almost all $x\in\Real^d$. Hence by Lemma \ref{difloc} $\mathrm{(a)}$ we have $g^{\Delta,q}_{\mathrm{loc}}f(x)<\infty$, for almost all $x\in\Real^d$. In fact, from the proof of Lemma \ref{difloc} it can be deduced that $\displaystyle\norm{t\partial_tP_t\left(\chi_N(x,\cdot)f(\cdot)\right)(x)}_{L_X^q((0,\infty),\frac{dt}{t})}<\infty$, for almost all $x\in\Real^d$. The arguments in the proof of  Theorem  \ref{herencia} can be used to conclude that $\norm{\sum_{k\geq1}\chi_{Q_k}(x)t\partial_tP_t\left(\chi_{2Q_k}f\right)(x)}_{L_X^q((0,\infty),\frac{dt}{t})}$ is finite for almost all $x\in\Real^d$. By the finite overlapping property of the balls $Q_k$ we get the finiteness almost every $x$ of each term $\norm{\chi_{Q_k}(x)t\partial_tP_t\left(\chi_{2Q_k}f
 \right)(x)}_{L_X^q((0,\infty),\frac{dt}{t})}$. On the other hand, observe that
\begin{align*}
    \norm{\chi_{Q_k}(x)t\partial_tP_t\left((1-\chi_{2Q_k})f\right)(x)}_{L_X^q((0,\infty),\frac{dt}{t})}^q &\leq C\int_0^\infty\left(\int_{\abs{x-y}>\rho(x_k)}\frac{t\norm{f(y)}_X}{(t+|x-y|)^{d+1}}~dy\right)^q\frac{dt}{t} \\
    &\leq C\norm{f}_{L^1_X(\Real^d)}^q\int_0^\infty\frac{t^q}{(t+\rho(x_k))^{(d+1)q}}~\frac{dt}{t} \\
    &\leq C_k\norm{f}_{L^1_X(\Real^d)}^q.
\end{align*}
Pasting together the last two thoughts we get that for every $k$ and almost every $x\in\Real^d$ the norm $\norm{\chi_{Q_k}(x)t\partial_tP_tf(x)}_{L_X^q((0,\infty),\frac{dt}{t})}$ is finite. Hence $\norm{t\partial_tP_tf(x)}_{L_X^q((0,\infty),\frac{dt}{t})}=g^{\Delta,q}f(x)$ is finite for almost all $x$.
\end{proof}



\end{document}